\documentclass[letterpaper, leqno]{article}
\usepackage{a4, amsxtra,amsmath,amsfonts,amscd, mathrsfs}
\usepackage[all]{xy}
\newcommand{\kkk}[1]{}

\newtheorem{thm}{Theorem}[section]
\newtheorem{prop}[thm]{Proposition}
\newtheorem{cor}[thm]{Corollary}
\newtheorem{lem}[thm]{Lemma}

\newtheorem{rem}[thm]{Remark}
\newtheorem{art}[thm]{}

\newcommand{\codim}{{\rm codim}}

\newcommand{\Div}{\rm div}

\newcommand{\Pic}{\rm Pic}
\newcommand{\ord}{\rm ord}
\newcommand{\Spec}{\rm Spec}
\newcommand{\Spf}{\rm Spf}

\newcommand{\ve}{\varepsilon}

\newcommand{\supp}{{\rm supp}}

\newcommand{\Acal}{{\mathscr A}}
\newcommand{\Ccal}{{\mathscr C}}

\newcommand{\Gcal}{{\mathscr G}}
\newcommand{\Hcal}{{\mathscr H}}

\newcommand{\Lcal}{{\mathscr L}}
\newcommand{\Mcal}{{\mathscr M}}
\newcommand{\Ncal}{{\mathscr N}}
\newcommand{\Ocal}{{\mathscr O}}

\newcommand{\Ucal}{{\mathscr U}}
\newcommand{\Xcal}{{\mathscr X}}

\newcommand{\qdop}{{\mathbb Q}}
\newcommand{\ndop}{{\mathbb N}}
\newcommand{\rdop}{{\mathbb R}}

\newcommand{\kdop}{{\mathbb K}}

\newcommand{\zdop}{{\mathbb Z}}

\newcommand{\metr}{\|\hspace{1ex}\|}

\newcommand{\g}{{\frak g}_{X,v}}
\newcommand{\gp}{{\frak g}_{X,v}^+}
\newcommand{\gh}{\hat{\frak g}_{X,v}}
\newcommand{\ghp}{\hat{\frak g}_{X,v}^+}

\newcommand{\proof}{\noindent {\bf Proof: \/}}

\newcommand{\qed}{{ \hfill $\square$}}

\newcommand{\val}{{\rm val}}

\newcommand{\ub}{{\mathbf u}}
\newcommand{\ubb}{{\overline{\mathbf u}}}

\newcommand{\mb}{{\mathbf m}}
\newcommand{\nb}{{\mathbf n}}

\newcommand{\xb}{{\mathbf x}}
\newcommand{\yb}{{\mathbf y}}

\newcommand{\Tor}{{\mathbb G}_m^n}
\newcommand{\rtor}{{\rdop^n/\Lambda}}

\newcommand{\Deltabar}{{\overline{\Delta}}}
\newcommand{\Ccalbar}{{\overline{\Ccal}}}

\newcommand{\valbar}{{\overline{\val}}}
\newcommand{\Xan}{{X_v^{\rm an}}}
\title{The Bogomolov conjecture \linebreak for totally degenerate abelian varieties}
\date{\today}
\author{Walter Gubler}
\begin{document}
\maketitle

\section{Introduction}

Let $K=k(B)$ be a function field of an integral projective variety $B$ over the algebraically closed field $k$ such that $B$ is regular in codimension $1$. The set of places $M_B$ is given by the prime divisors of $B$. We fix an ample class $\mathbf c$ on $B$. If we count every prime divisor $Y$ with weight $\deg_{\mathbf c}(Y)$, then the valuations $\ord_Y$ lead to a product formula on $K$ and hence to a theory of heights (see \cite{La} or \cite{BG}, Section 1.5). The algebraic closure of $K$ is denoted by $\overline{K}$.
For an abelian variety $A$ over $K$ which is totally degenerate at some place $v \in M_B$ (see \S 5 for definition), we will prove the {\it Bogomolov conjecture}:

\begin{thm} \label{Bogomolov conjecture} 
Let $X$ be a closed subvariety of $A$ defined over $\overline{K}$ which is not a translate of an abelian subvariety by a torsion point. For every ample symmetric line bundle $L$ on $A$, there is $\ve > 0$ such that
$$X(\ve):=\{P \in X(\overline{K}) \mid \hat{h}_L(P) \leq \ve\}$$
is not Zariski dense in $X$.
\end{thm}

Here, $\hat{h}_L$ is the N\'eron--Tate height with respect to $L$. For a number field $K$, Bogomolov \cite{Bo} conjectured this statement for a curve $X$ embedded in the Jacobian variety. This was proved by Ullmo \cite{Ul} and generalized by Zhang \cite{Zh2} to higher dimensional $X$ in an arbitrary abelian variety $A$. In the number field case, $A$ is not assumed to be totally degenerate at some finite place. 

For function fields however, the Bogomolov conjecture can't be true for arbitrary abelian varieties $A$ in the form given in Theorem \ref{Bogomolov conjecture}. Conjecturally, one has also to omit that $X=\varphi(Y_{\overline{K}})$ for a homomorphism $\varphi:C_{\overline{K}}\rightarrow A$ of an abelian variety $C$ over $k$ with a subvariety $Y$ also defined over $k$. It is very surprising that the function field case of Bogomolov's conjecture remained open. In fact, only some special cases for curves $X$ embedded in the Jacobian were proved by Moriwaki (\cite{Mo1}, \cite{Mo2}, \cite{Mo3}) and Yamaki (\cite{Ya}). Their proofs are based on analysis on the reduction graphs of a suitable semistable model of $X$.

The proof of Theorem \ref{Bogomolov conjecture} follows Zhang's proof replacing the complex analytic methods by tropical analytic geometry from \cite{Gu3} at the place $v$. Note that Moriwaki \cite{Mo4} proved the Bogomolov conjecture for finitely generated fields over $\qdop$ with respect to a set of (almost) absolute values which generalizes the number field situation but which is different from the classical function field case. Our arguments for Theorem \ref{Bogomolov conjecture} work also in Moriwaki's case (and hence for number fields) leading to a non-archimedean proof if $A$ is totally degenerate at a non-archimedean place $v$ (see Remark \ref{Moriwaki}).

The paper is organized as follows. In \S 2, we recall the theory of local heights and we introduce Chambert-Loir's measures. In \S 3, we resume global heights over function fields. The fundamental inequality and the relations to the successive minima are given in \S 4. They are completely analogous to the number field case. This is used in \S 5 to prove the tropical equidistribution theorem. The proof is more subtle than in the known equidistribution theorems as we have to vary the metric very carefully at the place $v$ (see Remark \ref{other proofs}). The Bogomolov conjecture is proved in \S 6 with the help of the tropical equidistribution theorem. Finally, we give some corollaries analogous to Zhang's paper \cite{Zh2}.

\vspace{3mm}

\centerline{\it Terminology}

In $A \subset B$, $A$ may be  equal to $B$. The complement of $A$ in $B$ is denoted by $B \setminus A$ \label{setminus} as we reserve $-$ for algebraic purposes. The zero is included in $\ndop$. The standard scalar product of $\ub,\ub' \in \rdop^n$ is denoted by  $\ub \cdot \ub':=u_1 u_1' + \dots + u_n u_n'$.

By a function field $K=k(B)$, we mean always the setting as at the beginning of the introduction with an ample class $\mathbf c$ fixed on $B$.
All occuring rings and algebras are commutative with $1$. If $A$ is such a ring, then the group of multiplicative units is denoted by $A^\times$. A variety over a field is a separated reduced scheme of finite type. For the degree of a map $f:X \rightarrow Y$ of irreducible varieties, we use either $\deg(f)$ or $[X:Y]$.  

\small
The author thanks J. Eckhoff, K. K\"unnemann and F. Oort for precious discussions, and the referee for his suggestions. Part of the research in this paper was done during a $5$ weeks stay at the CRM in Barcelona.
\normalsize

\section{Local heights}

In this section, $K$ denotes a field with a discrete valuation $v$. Let $\kdop_v$ be the completion of the algebraic closure of the completion of $K$. This is an algebraically closed field complete with respect to the absolute value $|\phantom{a}|_v$ extending the given absolute value $e^{-v}$ on $K$ (see \cite{BGR}, Proposition 3.4.1/3).

We will define local heights for a $d$-dimensional projective variety $X$ over $K$. For simplicity, we assume that $X$ is geometrically integral. In general, we can reduce to this special case by base change and linearity. First, we recall the facts needed from analytic and formal geometry. After a summary of local heights of varieties, we introduce Chambert-Loir's measures on $X$ and then we discuss canonical heights and canonical measures on abelian varieties. Most results hold more generally for complete varieties and for arbitrary valuations of height one (see \cite{Gu2} and \cite{Gu3}).

\begin{art} \rm \label{Berkovich analytic spaces}
Let $X_v$ be the base change of $X$ to $\kdop_v$ and let $X_v^{\rm an}$ be the {\it Berkovich analytic space} associated to $X_v$. Similarly, we proceed with line bundles on $X$ or morphisms. The GAGA theorems, well known in the archimedean situation, hold also in this situation. For more details about this functorial construction, we refer to \cite{Ber1}, 3.4.
\end{art}

\begin{art} \rm \label{admissible formal schemes}
The valuation ring of $\kdop_v$ is denoted by $\kdop_v^\circ$, $\kdop_v^{\circ \circ}:= \{\alpha \in \kdop_v \mid |\alpha|_v < 1\}$ and $\tilde{\kdop}_v:=\kdop_v^\circ/\kdop_v^{\circ \circ}$ is the residue field. 

A $\kdop^\circ$-algebra is called {\it admissible} if it is isomorphic to $\kdop^\circ\langle x_1, \dots , x_n \rangle / I$ for an ideal $I$ and if $A$ has no $\kdop^\circ$-torsion. An {\it admissible formal scheme} $\Xcal$ over $\kdop^\circ$ is a formal scheme which has a locally finite atlas of open subsets isomorphic to $\Spf(A)$ for  admissible $\kdop^\circ$-algebras $A$. The lack of $\kdop^\circ$-torsion is equivalent to flatness over $\kdop^\circ$ (see \cite{BL2}, \cite{BL3} for details). 

Then the {\it generic fibre} $\Xcal_v^{\rm an}$ of $\Xcal_v$ is the analytic space locally defined by the Berkovich spectrum of the affinoid algebra $A \otimes_{\kdop_v^\circ} \kdop_v$. Similarly, the {\it special fibre} $\tilde{\Xcal}_v$ of $\Xcal$ is the scheme over $\tilde{\kdop}_v$ locally defined by $\Spec(A/\kdop_v^{\circ\circ}A)$. By \cite{Ber2}, \S 1, there is a  surjective reduction map 
$$\Xcal_v^{\rm an} \rightarrow \tilde{\Xcal}_v, \, x \mapsto x(v).$$
\end{art}

\begin{art} \rm \label{models}
A {\it formal $\kdop_v^\circ$-model} of $X$ is an admissible formal scheme $\Xcal_v$ over $\kdop_v^\circ$ with $\Xcal_v^{\rm an}=X_v^{\rm an}$. For a line bundle $L$ on $X$, a formal $\kdop_v^\circ$-model of $L$ on $\Xcal_v$ is a line bundle $\Lcal_v$ on $\Xcal_v$ with generic fibre $\Lcal_v^{\rm an}=L_v^{\rm an}$.

In fact, $\Xcal_v$ is always dominated by the formal completion along the special fibre of a projective $\kdop_v^\circ$-model where every line bundle is ``difference'' of two very ample line bundles (see \cite{Gu2}, Proposition 10.5) and hence we may work algebraically. However, we will see that the analytic structure $X_v^{\rm an}$ is very important for the equidistribution and so it is natural to use formal models over $\kdop_v^\circ$ which allows more flexibility in choosing models.

A metric $\metr_v$ on $L_v^{\rm an}$ is said to be a {\it formal metric} if there is a formal $\kdop_v^\circ$-model $\Lcal_v$ of $L$ such that
for every formal trivialization $\Ucal_v$ of $\Lcal_v$ and every $s \in \Gamma(\Ucal_v,\Lcal_v)$ corresponding to $\gamma \in \Ocal_{\Xcal_v}(\Ucal_v)$, we have 
$\|s(x)\|=|\gamma(x)|$ on $\Ucal_v^{\rm an}$. 
The formal metric is called {\it semipositive} if there is $k \in \ndop \setminus \{0\}$ such that $ \Lcal_v^{\otimes k}$ is generated by global sections.
Note that this definition of semipositivity is more restrictive than the ones in \cite{Gu2}, \cite{Gu3}. It gives satisfactory results for projective varieties and it is enough to handle canonical local heights of symmetric ample line bundles on an abelian variety which is our main purpose. Moreover, it has the advantage that the corresponding local heights are non-negative.

A metric $\metr_v$ on $L_v^{\rm an}$ is called a {\it root of a formal metric} if there is a non-zero $k \in \ndop$ such that $\metr_v^{\otimes k}$ is a formal metric on $L^{\otimes k}$. 
\end{art}

\begin{art} \rm \label{isometry classes}
Let $\g$ be the set of isometry classes of formally metrized line bundles on $X$. It is the group with respect to the tensor product which is generated by the semigroup $\gp$ of semipositive isometry classes. We write the group $\g$ additively. 

We consider more generally metrics on a line bundle $L$ which are uniform limits of roots of semipositive formal metrics. Here, the distance of metrics on $L$ is given by 
$$d(\metr_v, \metr_v'):= \max_{x \in X_v^{\rm an}} | \log(\metr'_v/\metr_v)(x)|,$$
where $\metr'_v/\metr_v$ is evaluated at the section $1$ of $O_X$. The isometry classes of line bundles on $X$ endowed with such metrics form a semigroup $\ghp$. The {\it completion} $\gh$ of $\g$ is the group $\gh:=\ghp - \ghp$. 
\end{art}

\begin{art} \label{local height} \rm 
We consider pseudo-divisors $D_0, \dots , D_d$ on $X$ with $\g$-metrics $\metr_{j,v}$ on $O(D_j)_v^{\rm an}$ satisfying
\begin{equation} \label{suppcond}
\supp(D_0) \cap \dots \cap \supp(D_t)  = \emptyset.
\end{equation}
If the reader is not familiar with pseudo-divisors (see \cite{Fu}, \S 2), then the use of Cartier divisors is possible. This is more restrictive and formally more complicated, but it will be enough for our purposes. By intersection theory, we define a local height  of $X$. For simplicity, we explain it in the case of a projective $K^\circ$-model $\Xcal_v$ of $X$ with line bundles $\Lcal_0, \dots ,\Lcal_d$ inducing the metrics. Then the pseudo-divisor $D_j=(L_j,Y_j,s_j)$ gives rise to a pseudo-divisor ${\mathscr D}_j:=(\Lcal_j,Y_j \cup \tilde{\Xcal}_v,s_j)$ on $\Xcal_v$ and we define
the local height by
$$\lambda_{(D_0, \metr_{0,v}), \dots, (D_d, \metr_{d,v})}(X,v):= {\mathscr D}_0 \cdots {\mathscr D}_d \cdot \Xcal_v.$$ 
In general, we may use a suitable base change, but it is also possible to work with refined intersection theory of pseudo-divisors on formal $\kdop_v^\circ$-models (see \cite{Gu2}, \S 5). The extension of local heights to $\gh$-metrics is handled in:
\end{art}

\begin{thm} \label{properties of local height}
Let $D_0, \dots , D_d$ be pseudo-divisors on $X$ with \eqref{suppcond} and let $\metr_{0,v}, \dots , \metr_{d,v}$ be $\gh$-metrics on $O(D_0), \dots, O(D_d)$. Then there is a unique  $\lambda_{(D_0, \metr_{0,v}), \dots, (D_d, \metr_{d,v})}(X,v) \in \rdop$ with the following properties: 
\begin{itemize}
\item[(a)] It is multilinear and symmetric in $(D_0, \metr_{0,v}), \dots, (D_d, \metr_{d,v})$.
\item[(b)] Let  $\varphi:X' \to X$ be a morphism of geometrically integral $d$-dimensional varieties, then 
\begin{multline*}
\lambda_{(\varphi^*(D_0),\varphi^*\metr_{0,v}), \dots, (\varphi^*(D_d),\varphi^*\metr_{d,v})}(X',v) \\
= \deg(\varphi) \lambda_{(D_0, \metr_{0,v}), \dots, (D_d, \metr_{d,v})}(X,v).
\end{multline*}
\item[(c)] If $D_0=\Div(f)$ for a rational function $f$ on $X$ and if $\metr_{v,0}$ is the trivial metric on $O(D_0)=O_X$, then
$$\lambda_{(D_0, \metr_{0,v}), \dots, (D_d, \metr_{d,v})}(X,v) = \log|f(Y)|_v,$$ 
where $Y$ is a representative of $D_1 \dots D_t.X \in CH_0 \left(\supp(D_1) \cap \dots \cap \supp(D_t)  \right)$ and where the right hand side is defined by linearity in the components of $Y_{\kdop_v}$.
\item[(d)] Suppose that $\metr_{1,v}, \dots ,\metr_{d,v}$ are all $\ghp$-metrics and that $\metr_{0,v}'$ is a second $\gh$-metric on $O(D_0)$. Then
\begin{multline*}
|\lambda_{(D_0, \metr_{0,v}), \dots, (D_d, \metr_{d,v})}(X,v) - 
\lambda_{(D_0, \metr_{0,v}'), \dots, (D_d, \metr_{d,v})}(X,v) | \\
\leq d(\metr_{0,v}, \metr_{0,v}') \deg_{D_1, \dots , D_d}(X). 
\end{multline*}
\item[(e)] If $\metr_{0,v}, \dots , \metr_{d,v}$ are all $\ghp$-metrics, then the local height is non-negative.
\item[(f)] If $\metr_{0,v}, \dots , \metr_{d,v}$ are formal metrics, then the local height is given by \ref{local height}.
\end{itemize}
\end{thm}

For a proof, we refer to \cite{Gu2}, Theorem 10.6.

\begin{art} \rm \label{Chern forms}
A continuous function $g$ on $X_v^{\rm an}$ induces a metric $\metr_{g,v}$ on $O_X$ determined by $||1||_{g,v}:=e^{-g}$. We assume that $||1||_{g,v}$ is a root of a formal metric. By \cite{Gu1}, Theorem 7.12, such functions form a dense $\qdop$-subspace of $C(X_v^{\rm an})$. 

For $\gh$-metrized line bundles $(L_1,\metr_{1,v}), \dots , (L_d,\metr_{d,v})$ on $X$, we get metrized pseudo-divisors $\hat{D}_j:= (L_j, \emptyset, 1, \metr_{j,v})$ and $\hat{O}^g:=(O_X,  \emptyset,1, \metr_g)$.  Then we define
$$\int_{X^{\rm an}} g \,c_1(L_1,\metr_{1,v}) \wedge \dots \wedge c_1(L_d,\metr_{d,v})
:= \lambda_{\hat{O}^g, \hat{D}_1, \dots, \hat{D}_d}(X,v).$$ 
By Theorem \ref{properties of local height}, this extends uniquely to a continuous functional on $C(X_v^{\rm an})$ and hence we get a regular Borel measure $c_1(L_1,\metr_{1,v}) \wedge \dots \wedge c_1(L_d,\metr_{d,v})$ on $X_v^{\rm an}$. 
These measures were first introduced by Chambert-Loir (see \cite{Ch}, \S2). For the following result, which follows from Theorem \ref{properties of local height}, we refer to \cite{Gu3}, Proposition 3.9, Corollary 3.11 and Proposition 3.12.
\end{art}

\begin{cor} \label{Chern properties}
The regular Borel measure $c_1(L_1,\metr_{1,v}) \wedge \dots \wedge c_1(L_d,\metr_{d,v})$ has the following properties:
\begin{itemize}
\item[(a)] It is multilinear and symmetric in  $(L_1,\metr_{1,v}), \dots , (L_d,\metr_{d,v})$.
\item[(b)] Let $\varphi:X' \rightarrow X$ be a morphism of geometrically integral projective varieties of dimension $d$, then
$$\varphi_* \left( c_1(\varphi^*L_1,\varphi^*\metr_{1,v}) \wedge \dots \wedge c_1(\varphi^*L_d,\varphi^*\metr_{d,v})\right) 
= \deg(\varphi) c_1(L_1,\metr_{1,v}) \wedge \dots \wedge c_1(L_d,\metr_{d,v}).$$
\item[(c)] If all the metrics are semipositive, then $c_1(L_1,\metr_{1,v}) \wedge \dots \wedge c_1(L_d,\metr_{d,v})$ is a positive measure of total measure $\deg_{L_1, \dots , L_d}(X)$ which depends continuously on the metrics with respect to the weak topology on the set of positive regular Borel measures on $\Xan$.
\item[(d)] If all metrics are formal, induced by formal $\kdop_v^\circ$-models $\Lcal_1, \dots, \Lcal_d$ on the formal $\kdop_v^\circ$-model $\Xcal$ of $X$ with reduced special fibre $\tilde{\Xcal}$, then
$$c_1(L_1,\metr_{1,v}) \wedge \dots \wedge c_1(L_d,\metr_{d,v})=\sum_Y \deg_{\tilde{\Lcal}_1, \dots, \tilde{\Lcal}_d}(Y) \delta_{\xi_Y},$$
where $Y$ ranges over the irreducible components of $\tilde{\Xcal}$ and $\delta_{\xi_Y}$ is the Dirac measure in the unique point $\xi_Y \in X_v^{\rm an}$ with reduction $\xi_Y(v)$ equal to the generic point of $Y$.
\end{itemize}
\end{cor}

\begin{rem} \rm \label{multiple of metric}
If we replace the metric on $L_j$ by a $|\kdop_v^\times|_v$-multiple, then $c_1(L_1,\metr_{1,v}) \wedge \dots \wedge c_1(L_d,\metr_{d,v})$ remains the same. Indeed, we may check this for a function $g$ associated to a root of a formal metric on $O_X$ and then the claim follows from Theorem \ref{properties of local height}(a) and (c).
\end{rem}

\begin{art} \rm \label{canoncial local heights}
Let $A$ be an abelian variety over $K$. We consider an ample symmetric rigidified line bundle $(L,\rho)$ on $A$, where $\rho \in L(K) \setminus \{0\}$. Then there exists a canonical ${\hat{\frak g}_{A,v}^+}$-metric $\metr_{\rho,v}$ on $L_v^{\rm an}$ which behaves well with respect to tensor product and homomorphic pull-back (see \cite{BG}, Theorem 9.5.7 and its proof). 

Let $X$ be a closed subvariety of dimension $d$. For $j=0, \dots, d$, let $D_j$ be a pseudo-divisor on $A$ and let $\rho_j$ be a rigidification on $O(D_j)$. If $\supp(D_0) \cap \dots \cap \supp(D_d) \cap X = \emptyset$, then  we get a {\it canonical local height}
$$\hat{\lambda}_{(D_0,\rho_0), \dots, (D_d,\rho_d)}(X,v):=
\lambda_{(D_0,\metr_{\rho_0,v}), \dots, (D_d,\metr_{\rho_d,v})}(X,v).$$ 
Rigidified ample symmetric line bundles $(L_1,\rho_1), \dots, (L_d,\rho_d)$ induce a {\it canonical measure}
$$c_1(L_1|_X,\metr_{\rho_1,v}) \wedge \dots \wedge c_1(L_d|_X,\metr_{\rho_d,v})$$
on $X_v^{\rm an}$. By Remark \ref{multiple of metric}, the canoncial measure does not depend on the choice of the rigidifications. We leave it to the reader to transfer Theorem \ref{properties of local height} and Corollary \ref{Chern properties} to  canonical local heights and measures.
\end{art}

\section{Global heights}

In this section, we resume the theory of global heights over the function field $K=k(B)$. For simplicity, we restrict to the case of an irreducible projective variety $X$ over $K$ of dimension $d$. 

\begin{art} \rm \label{model heights}
Recall from the introduction that we have fixed an ample class $\mathbf c$ on $B$ such that  the prime divisors $Y$ on $B$ are weighted by $\deg_{\mathbf c}(Y)$. More generally, $\deg_{\mathbf c}(Z)$ denotes the degree of a cycle $Z$ on $B$ with respect to $\mathbf c$. For  $Z$ of pure dimension $e$, it is equal to the intersection number $\mathbf c^{\cdot e}\cdot Z$ and hence it depends only on the rational equivalence class of $Z$. In particular, we may define $\deg_{\mathbf c}({\cal B}):=\deg_{\mathbf c}(c_1({\cal B}))$ for a line bundle ${\cal B}$ on $B$.

Let $\Xcal$  be an irreducible projective variety over $B$ with generic fibre $X$ and let $\Lcal$ be a line bundle on $\Xcal$ with generic fibre $L$ on $X$. Then the {\it height} of  $P \in X(\overline K)$ with respect to $\Lcal$ is given by 
$$h_{\Lcal}(P):=\frac{1}{[K(P):K]} \deg_{\mathbf c'}(\varphi_P^* \Lcal).$$
Here, $P$ is rational  over some finite extension $F/K$ which is a function field  $F=k(B')$ for the normalization $B'$ of $B$ in $F$. Moreover, $\mathbf c'$ denotes the pull-back of $\mathbf c$ to $B'$. Note that the canonical rational map $\varphi_P:B' \dashrightarrow X$  induced by $P$ is defined in codimension $1$. To stress the similarity to arithmetic intersection theory, let  
$$h_{\Lcal}(X):= \deg_{\mathbf c}\left(\pi_*\left(c_1(\Lcal)^{d+1}.\Xcal\right)\right)$$
where $\pi:\Xcal \rightarrow B$ is the morphism of structure. More generally, the {\it height} of $X$ with respect to line bundles $\Lcal_0, \dots, \Lcal_d$ is defined by
$$h_{\Lcal_0, \dots, \Lcal_d}(X):= \deg_{\mathbf c}\left(\pi_*\left(c_1(\Lcal_0) \dots c_1(\Lcal_d).\Xcal\right)\right).$$
\end{art}

\begin{art} \rm \label{M-metric}
To define canonical global heights of subvarieties of an abelian variety over $K$, it is necessary to allow the local heights from Section 2 at finitely many places. For this, we define an {\it admissible $M_B$-metric} $\metr$ on a line bundle $L$ of $X$ as a family $\metr:= (\metr_v)_{v \in M_B}$ of metrics $\metr_v$ on $L_v^{\rm an}$ with the following properties:
\begin{itemize}
\item[(a)] There is an open dense subset $V$ of $B$, an irreducible projective variety $\Xcal$ over $B \setminus V$ with generic fibre $X$, a non-zero $N \in \ndop$ and a line bundle $\Lcal$ on $\Xcal$ such that $L^{\otimes N} = \Lcal |_X$ and such that $\metr_v^{\otimes N}$ is the formal metric associated to $\Lcal$ at all places $v \in M_B \cap V$.
\item[(b)] For $v \in M_B \setminus V$, $\metr_v$ is a $\gh$-metric on $L_v^{\rm an}$.
\end{itemize}
The admissible $M_B$-metric $\metr$ on $L$ is called {\it semipositive} if we can choose $\Lcal$ in (a) to be generated by global sections and if $\metr_v \in \ghp$ for all $v \in S$.
\end{art}

\begin{art} \rm \label{global M-height}
Let $\overline{L}_0 ,\dots ,\overline{L}_d$ be line bundles on $X$ endowed with admissible $M_B$-metrics. For $j=0, \dots, d$, we choose a pseudo-divisor $D_j$ with $L_j=O(D_j)$ and with
$$\supp(D_0) \cap \dots \cap \supp(D_d) = \emptyset.$$
Then we define the height of $X$ with respect to $\overline{L}_0 ,\dots ,\overline{L}_d$ by
$$h_{\overline{L}_0 ,\dots ,\overline{L}_d}(X):=
\sum_{v \in M_B} \lambda_{(D_0, \metr_{0,v}), \dots , (D_d, \metr_{d,v})}(X,v).$$ 
By (a), only finitely many $v \in M_B$ give a non-zero contribution. Moreover, Theorem \ref{properties of local height}(c) and the product formula show that the height is independent of the choice of $D_0, \dots ,D_d$. 
\end{art}

\begin{thm} \label{properties of global heights}
The height $h_{\overline{L}_0 ,\dots ,\overline{L}_d}(X)$ is uniquely determined by the following properties:
\begin{itemize}
\item[(a)] It is multilinear and symmetric in $\overline{L}_0 ,\dots ,\overline{L}_d$.
\item[(b)] If $\varphi:X' \rightarrow X$ is a morphism of irreducible $d$-dimensional projective varieties, then
$$h_{\varphi^*\overline{L}_0, \dots, \varphi^*\overline{L}_d}(X')=
\deg(\varphi) h_{\overline{L}_0 ,\dots ,\overline{L}_d}(X).$$
\item[(c)] If the $M_B$-metrics of $\overline{L}_1 ,\dots ,\overline{L}_d$ are semipositive and if we consider two $M_B$-admissible metrics on $L_0$ inducing the metrized line bundles $\overline{L}_0$ and $\overline{L}_0'$, then
$$h_{\overline{L}_0',\overline{L}_1 ,\dots ,\overline{L}_d}(X)
- h_{\overline{L}_0 ,\dots ,\overline{L}_d}(X)=O(\deg_{L_1, \dots,L_d}(X)).$$
\item[(d)] If the $M_B$-metrics of $
\overline{L}_0 ,\dots ,\overline{L}_d$ are semipositive, then
$$h_{\overline{L}_0 ,\dots ,\overline{L}_d}(X)\geq 0.$$
\item[(e)] If the $M_B$-metrics of $\overline{L}_0 ,\dots ,\overline{L}_d$ are induced by line bundles $\Lcal_0, \dots, \Lcal_d$ as in \ref{model heights}, then
$$h_{\overline{L}_0 ,\dots ,\overline{L}_d}(X)=h_{\Lcal_0, \dots, \Lcal_d}(X).$$
\end{itemize}
\end{thm}

\proof See \cite{Gu2}, 11.7. \qed

\begin{art} \label{Neron-Tate height} \rm 
Let $A$ be an abelian variety over $K$. For a rigidified ample symmetric line bundle $(L,\rho)$ on $A$, the canonical metrics $\metr_{\rho,v}$ form an admissible $M_B$-metric $\metr_\rho$ which we call {\it canonical}. Indeed, $A$ extends to a projective group scheme $\Acal$ over an open dense subset $V$ of $B$ and we may assume that $L$ extends to a symmetric line bundle $\Lcal$ on $\Acal$ which is relatively ample over $V$. By shrinking $V$ further, we may assume that some tensor power of $\Lcal$ will be generated by global sections and hence $\metr_\rho$ will be semipositive.

Let $L_0, \dots, L_d$ be ample symmetric line bundles on $A$. We endow them with canonical metrics $\metr_{\rho_j,v}$ as above. The {\it N\'eron--Tate height} of an irreducible $d$-dimensional closed subvariety $X$ of $A$ with respect to $L_0, \dots,L_d$ is defined by 
$$\hat{h}_{L_0,\dots,L_d}(X):= h_{(L_0,\metr_{\rho_0,v}), \dots , (L_d,\metr_{\rho_d,v})} (X).$$
By the product formula, this does not depend on the choice of the rigidification $\rho_j$. The case $d=0$ yields a N\'eron--Tate height $\hat{h}_L(P)$ for $P \in A(\overline{K})$ similarly as in \ref{model heights}.
\end{art}

\begin{thm} \label{characterization Neron-Tate}
The N\'eron--Tate height $\hat{h}_{L_0, \dots, L_d}(X)$ is non-negative and it is uniquely characterized by the following properties:
\begin{itemize}
\item[(a)] It is multilinear and symmetric in the ample symmetric line bundles $L_0, \dots, L_d$.
\item[(b)] If $\varphi:A' \rightarrow A$ is a homomorphism of abelian varieties and if $X'$ is an irreducible $d$-dimensional closed subvariety of $A'$ mapping onto $X$, then
$$\hat{h}_{\varphi^*{L}_0, \dots, \varphi^*{L}_d}(X')=
\deg(\varphi) \hat{h}_{{L}_0 ,\dots ,{L}_d}(X).$$
\item[(c)] For $m \in \zdop$, we have $m^{2d+2} \hat{h}_{{L}_0 ,\dots ,{L}_d}(X)=[X:mX] \hat{h}_{{L}_0 ,\dots ,{L}_d}(mX)$.
\item[(d)]  If the admissible $M_B$-metrics of $\overline{L}_0 ,\dots ,\overline{L}_d$ are semipositive, then
$${h}_{\overline{L}_0 ,\dots ,\overline{L}_d}(X)-
\hat{h}_{L_0, \dots ,L_d}(X)
=\sum_{j=0}^d O(\deg_{L_0, \dots,L_{j-1},L_{j+1}, \dots, L_d}(X)).$$
\end{itemize}
\end{thm}

\proof See \cite{Gu2}, Theorem 11.18. \qed

\begin{rem} \rm \label{Tate approach}
If $\metr_0, \dots, \metr_d$ are any admissible semipositive $M_B$-metrics on $L_0,\dots,L_d$, then it follows that 
$$\hat{h}_{{L}_0 ,\dots ,{L}_d}(X)= \lim_{m \to \infty} m^{-2(d+1)} \hat{h}_{(L_0,\metr_0), \dots ,(L_d,\metr_d)}(X).$$
This is the higher dimensional analogue of Tate's approach to N\'eron--Tate heights.
\end{rem} 

\begin{rem} \rm \label{generalizations of global heights}
The results of this section (except the explicit description in \ref{model heights} may be generalized to complete varieties over arbitrary fields with product formula. Moreover, the N\'eron--Tate height may be defined with respect to  arbitrary line bundles on $A$. For details, we refer to \cite{Gu2}, \S 11.
\end{rem}

\section{The fundamental inequality}

Let $K=k(B)$ be a function field. In this section, we relate the height of an irreducible projective $d$-dimensional variety $X$ over $K$ to its successive minima. In the number field case, these results are due to Zhang \cite{Zh1}.  

Let $\overline{L}$ be an ample line bundle on $X$ endowed with an admissible semipositive $M_B$-metric $\metr= (\metr_v)_{v \in M_B}$. For $i=1, \dots, d+1$, we consider the {\it successive minima}
$$e_i(X,\overline{L}) := \sup_Y \inf_{P \in X(\overline{K}) \setminus Y} h_{\overline{L}}(P),$$
where $Y$ ranges over all closed subsets of codimension $i$ in $X$.

\begin{lem}[fundamental inequality] \label{fundamental inequality} Under the hypothesis above, we have
$$\frac{h_{\overline L}(X)}{(d+1)\deg_L(X)} \leq e_1(X,\overline L).$$
\end{lem}

\proof To prove the fundamental inequality, we may assume that $B$ is a smooth projective curve. Indeed, we may first assume that $\mathbf c$ is very ample by passing to a tensor power and then we may replace $B$ by the generic curve obtained by intersecting $B$ with generic hyperplane sections of $\mathbf c$. This process does not change the relevant quantities (see \cite{BG}, 14.5.8).

By Theorem \ref{properties of local height}(d), both sides are continuous with respect to $\ghp$-metrics on $L$ at a given place $v \in B$. Moreover, it is enough to prove the claim for a suitable tensor power of $\overline{L}$. By definition of admissible semipositive metrics and after a suitable finite base change, it is enough to prove the claim for a formal $B$-metric $\metr$ on $L$ induced by a line bundle $\Lcal$ on an irreducible projective variety $\Xcal$ over $B$ with generic fibre $X$ such that $\Lcal$ is locally over $B$ generated by global sections and with $\Lcal|_X=L$. 

For an ample line bundle $\Hcal$ on $\Xcal$ and non-zero $N \in \ndop$, $\Lcal_N:=\Lcal^{\otimes N} \otimes \Hcal$ is relatively ample. If the fundamental inequality holds for all $\Lcal_N$, then we deduce, by $N \to \infty$ and using multilinearity, the fundamental inequality for $\Lcal$. 

Hence we may assume that $\Lcal$ is relatively ample over $B$. Using Snapper polynomials as in \cite{Fu}, Example 18.3.6, we conclude that the Euler characteristic $\chi(\Xcal, \Lcal^{\otimes n})$ is a numerical polynomial in $n \in \zdop$ with leading term
\begin{equation} \label{leading term}
\chi(\Xcal, \Lcal^{\otimes n}) = \frac{h_{\Lcal}(X)}{(d+1)!} n^{d+1} + O(n^d).
\end{equation}
 Under the above hypothesis, we prove the following Riemann--Roch result:

\begin{lem} \label{Riemann--Roch}
If $h_{\Lcal}(X) > 0$, then $H^0(\Xcal,\Lcal^{\otimes n}) \neq {0}$ for some sufficiently large $n \in \ndop$.
\end{lem}

\proof Since $\Lcal$ is relatively ample, we have 
$$H^i(\Xcal_t, \Lcal^{\otimes n})={0}$$
for all $t \in B$, $i>0$ and $n$ sufficiently large. By \cite{Mu2}, II.5, or \cite{Mi}, Theorem 4.2, the higher direct images with respect to the morphism of structure $\pi:\Xcal \rightarrow B$ satisfy
$$R^i\pi_*( \Lcal^{\otimes n})={0}$$
for every $i >0$ and as always $n$ sufficiently large. By \cite{Ha}, Exercise III.8.1, we get
$$H^i(\Xcal, \Lcal^{\otimes n})=H^i(B,\pi_* \Lcal^{\otimes n})$$
for every $i \geq 0$. As $B$ is a curve, the cohomology groups vanish for $i \geq 2$ and hence
$$\chi(\Xcal,\Lcal^{\otimes n}) \leq \dim H^0(\Xcal, \Lcal^{\otimes n}).$$
Now Lemma \ref{Riemann--Roch} is a consequence of \eqref{leading term}. \qed

{\it Continuation of the proof of the fundamental inequality:}

We may replace $\mathbf c$ by algebraically equivalent classes and tensor powers, hence we may assume that $\deg(\mathbf c)=1$. This is usually the standard hypothesis for function fields of curves. We denote by $\Mcal$ a line bundle on $B$ with isomorphism class $\mathbf c$. For $a \in \qdop$, we set
$$\Lcal(a):= \Lcal \otimes \pi^*(\Mcal)^{\otimes a} \in \Pic(X) \otimes_\zdop \qdop.$$
Then $\Lcal(a)$ remains relatively ample. We choose 
$$a:= -\frac{h_\Lcal(X)}{(d+1)\deg_L(X)} + \ve$$
for any positive rational $\ve$. Then projection formula and $\deg(\Mcal)=1$ show
$$\deg_{\Lcal(a)}(\Xcal)=c_1(\Lcal(a))^{d+1} \cdot \Xcal
=h_\Lcal(X) + (d+1) a  \deg_\Lcal(X) >0.$$
By Lemma \ref{Riemann--Roch}, there is a non-trivial global section $s$ of $\Lcal(a)^{\otimes n}$ for some sufficiently large $n \in \ndop$. Since $\Lcal(a)^{\otimes n}$ is a line bundle, the restriction of $s$ to the generic fibre $X$ does not vanish identically. We choose $P \in X(\overline K)$ outside the support of $\Div(s)$. Then 
$$h_\Lcal(P)+a = h_{\Lcal(a)}(P) = \frac{1}{[K(P):K]} \deg \Div(\varphi_P^*(s)) \geq 0.$$
We conclude that
$$h_\Lcal(P) \geq -a= \frac{h_\Lcal(X)}{(d+1)\deg_L(X)} - \ve.$$
As we may choose $\ve$ arbitrarily small, we get the fundamental inequality. \qed

Similarly as in Section 5 of \cite{Zh1}, there is also a converse to the fundamental inequality:

\begin{prop} \label{converse}
Under the same hypothesis as in the fundamental inequality, we have
$$e_1(X,\overline{L}) + \dots + e_{d+1}(X,\overline{L}) \leq \frac{h_{\overline{L}}(X)}{\deg_L(X)}.$$
\end{prop}

\proof Again, we may assume that $B$ is a smooth projective curve and that the metric is induced by a relatively ample $\Lcal$. 
Similarly as in the proof of the fundamental inequality, we consider the line bundle 
$$\Ncal := \Lcal(-e_{d+1}(X,\Lcal)+\ve)$$
for some $\ve > 0$. By construction, the line bundle $\Lcal(-e_{d+1}(X,\Lcal))$ is numerically effective  and by the Nakai--Moishezon numerical criterion (\cite{Kl}, Theorem III.1), we conclude that $\Ncal$ is ample. Then the proof of the second inequality of Theorem 5.2 in \cite{Zh1} translates to the function field case immediately. \qed

\begin{cor} \label{height zero} 
We have $h_{\overline L}(X)=0$ if and only if $e_1(X,\overline L)=0$.
\end{cor}

\section{The tropical equidistribution theorem}

We first recall the basic facts from tropical analytic geometry for closed subvarieties $X$ of a totally degenerate abelian variety $A$ (see \cite{Gu3} for more details). Then we prove the equidistribution theorem for the tropical variety associated to $X$. 

We fix a place $v \in M_B$ of the function field $K=k(B)$. 
Let $A$ be an abelian variety over $K$ which is totally degenerate at the place $v$, i.e. $A_v^{\rm an}= (\Tor)_v^{\rm an}/M$ for a lattice $M$ in $\Tor(\kdop_v)$. The latter is a discrete subgroup of $T:=(\Tor)_v^{\rm an}$ which maps isomorphically onto a complete lattice $\Lambda$ by the map
$$\val:T \longrightarrow \rdop^n, \quad p \mapsto 
( v(x_1), \dots, v(x_n) ).$$
Here, we always fix coordinates $x_1, \dots, x_n$ of $\Tor$. Note that $\val$ is continuous.

\begin{art} \rm \label{Mumford's construction} 
Let $\pi:\rdop^n \rightarrow \rtor, \, \ub \mapsto \ubb$, be the quotient map. The usual definitions from convex geometry are translated to $\rtor$, e.g. a polytope $\Deltabar$ of $\rtor$ is given by a polytope $\Delta$ of $\rdop^n$ such that $\pi$ maps $\Delta$ bijectively onto $\Deltabar$. The map $\val$ descends to a continuous map $\valbar: A_v^{\rm an} \rightarrow \rtor$. Note that the valuation $v$ is discrete and hence the lattice $\Lambda$ is defined over $\qdop$. Let $\Ccalbar$ be a rational polytopal decomposition $\Ccalbar$ of $\rtor$, where rational means that the polytopes are defined over $\qdop$. 
By Mumford's construction (see \cite{Mu1}, \S 6, or \cite{Gu3}, \S 7), there is a canonical formal $\kdop_v^\circ$-model $\Acal_v$ of $A$ associated to $\Ccalbar$.
\end{art}

\begin{art} \rm \label{toric line bundles}
The line bundles on $A_v^{\rm an} = (\Tor)^{\rm an}_{\kdop_v}/M$ may be described in the following way (see \cite{FvdP}, Ch. VI, and \cite{BL1}, \S 2, for details). 

Let $L$ be a line bundle on $A_v^{\rm an}$. The pull-back to $T$ with respect to the quotient morphism $p$ is trivial and will be identified with $T \times {\kdop_v}$. It is given by a cocycle $\gamma \mapsto Z_\gamma$ of $H^1(M, \Ocal(T)^\times)$ and $L= (T \times {\kdop_v})/M$ where the quotient is with respect to the $M$-action
$$M \times \left( T \times {\kdop_v} \right)  \longrightarrow 
T \times {\kdop_v}, \quad (\gamma, (\xb, \alpha)) \mapsto (\gamma \cdot \xb, Z_\gamma(\xb)^{-1} \alpha).$$
The cocycle has the form $Z_\gamma(\xb)=d_\gamma \cdot \sigma_\gamma(\xb)$, where $\gamma \mapsto \sigma_\gamma$ is a homomorphism of $M$ to the character group $\check{T}$ and where $d_\gamma \in {\kdop_v}^\times$ satisfies
\begin{equation} \label{quadratic}
d_{\gamma \rho} \cdot d_\gamma^{-1} \cdot d_\rho^{-1} = \sigma_\rho(\gamma) \quad (\gamma, \rho \in M). 
\end{equation}
By the isomorphism $M \stackrel{\val}{\rightarrow} \Lambda$, we get a unique symmetric bilinear form $b$ on $\Lambda$ characterized by
$$b(\val(\gamma),\val(\rho)) = v\left( \sigma_\rho(\gamma) \right).$$
Then $b$ is positive definite on $\Lambda$ if and only if $L$ is ample. 
Note that the cocycle $Z_\gamma$ factors over $\rdop^n$, i.e. for every $\lambda = \val(\gamma) \in \Lambda$, there is a unique real function $z_\lambda$ on $\rdop^n$ such that 
$$z_{\lambda}(\val(\xb))= v (Z_\gamma(\xb)) \quad (\gamma \in M, \, \xb \in T).$$
Let $\Acal_v$ be the formal ${\kdop_v}^\circ$-model of $A$ associated to a given rational polytopal decomposition $\Ccalbar$ of $\rtor$. Let $f$ be a continuous function on $\rdop^n$ such that $f|_\Delta$ is an affine function over $\qdop$ of every $\Delta \in \Ccal$ and with 
\begin{equation} \label{periodicity}
f(\ub + \lambda) = f(\ub) + z_\lambda(\ub) \quad (\lambda \in \Lambda, \ub \in \rdop^n).
\end{equation}
Then $f$ induces a unique metric $\metr$ on $L$  with $f \circ \val = v \circ p^* \|1\|$ on $T$ (see \cite{Gu3}, Proposition 6.6). Moreover, $\metr$ is a root of a formal line bundle $\Lcal$ on $\Acal_v$ which has  ample reduction if and only if $f$ is a strongly polyhedral convex function with respect to $\Ccal$ (\cite{Gu3}, Corollary 6.7). Here, a convex function $f$ is called {\it strongly polyhedral with respect to $\Ccal$} if the $n$-dimensional $\Delta \in \Ccal$ are the maximal subsets of $\rdop^n$ where $f$ is affine.  
\end{art}

\begin{lem} \rm \label{barycentric model}
Let $L$ be an ample line bundle on $A$. Then there is a triangulation $\Ccalbar$ of $\rtor$ such that for all $j \in \ndop$, we have:
\begin{itemize}
\item[(a)] $\Ccal_j := 2^{-j} \Ccal$ is subdivided by $\Ccal_{j+1}$;
\item[(b)] the triangulation $\Ccalbar_j$ of $\rtor$ induces a formal $\kdop_v^\circ$-model $\Acal_v^{(j)}$ of $A$;
\item[(c)] there is $m_j \in \ndop$ and a formal $\kdop_v^\circ$-model $\Lcal_v$ of $L^{\otimes m_j}$ on $\Acal_v^{(j)}$ such that $\tilde{\Lcal}_v$ is ample.
\end{itemize}
\end{lem}

\proof By linear algebra, $\Lambda$ contains linearly independent vectors $b_1, \dots, b_n$ which are orthogonal with respect to the positive definite bilinear form $b$ associated to $L$. There is a non-zero $N \in \ndop$ such that
$$\Lambda \subset \Lambda':=\zdop b_1' + \dots + \zdop b_n', \quad b_j':=\frac{1}{N}b_j.$$
Then the fundamental lattice $F:=[0,1]b_1'+ \dots + [0,1]b_n'$ of $\Lambda'$ is a cuboid with respect to $b$. The {\it complete barycentric subdivision} of the $\Lambda'$-translates of $F$ leads to a triangulation $\Ccal$ of $\rdop^n$ (see \cite{Le}, 14.5.4). The simplex $\Delta$ with vertices
\begin{equation} \label{barycentric}
{\mathbf 0}, \quad \frac{1}{2}b_1', \quad \frac{1}{2}b_1'  + \frac{1}{2}b_2', \quad \dots ,\quad \frac{1}{2}b_1' + \dots + \frac{1}{2}b_n'
\end{equation}
is a typical maximal member of  $\Ccal$ and every other maximal simplex of $\Ccal$ is obtained from $\Delta$ by the following operations:
\begin{itemize}
 \item[(I)] replacing $b_1', \dots ,b_n'$ by $\pm b_1', \dots , \pm b_n'$;
 \item[(II)] permuting $b_1', \dots , b_n'$;
 \item[(III)] translating $\Delta$ by a vector $\lambda' \in \Lambda'$.
\end{itemize}
Clearly, $\Ccal_j$ induces a triangulation $\Ccalbar_j$ of $\rtor$ and Mumford's construction yields (b). It is an elementary exercise to check (a) for the barycentric subdivision of a cuboid. It is enough to prove (c) for $j=0$. We have to find  a rational strongly polyhedral convex function $f$ with respect to $\Ccal$ satisfying \eqref{periodicity}. Note that 
\begin{equation} \label{cycle and b}
z_\lambda(\ub)=z_\lambda({\mathbf 0})+b(\lambda, \ub)
\end{equation}
for $\lambda \in \Lambda$, $\ub \in \rdop^n$. By \eqref{quadratic}, $z_\lambda({\mathbf 0})$ is a quadratic function on $\Lambda$ and hence $z_\lambda({\mathbf 0})=q(\lambda)+\ell(\lambda)$ for the quadratic form $q(\lambda)=\frac{1}{2}b(\lambda,\lambda)$ and a linear form $\ell$ on $\Lambda$ (see \cite{BG}, 8.6.5). We extend $q$ and $\ell$ to $\rdop^n$ to get a strictly convex function $q+\ell$ on $\rdop^n$ which satisfies \eqref{periodicity} as a consequence of \eqref{cycle and b}.  

In a first try,  let $f$ be the function on $\rdop^n$ which is affine on every maximal simplex $\Delta$ of $\Ccal$ and agrees with $q+\ell$ on the vertices of $\Delta$. Obviously, $f$ also satisfies \eqref{periodicity}. We have to check that $f$ is a strongly polyhedral convex function with respect to $\Ccal$. Let $\Delta$ and $\Delta'$ be maximal simplices of $\Ccal$ such that $\sigma := \Delta \cap \Delta'$ is a simplex of codimension $1$ in $\Delta$ (and $\Delta'$). Let $\ub_0$ (resp. $\ub_0'$) be the remaining vertex of $\Delta$ (resp. $\Delta'$) outside of $\sigma$. Since $f$ is piecewise affine with respect to $\Ccal$, it is enough to show 
\begin{equation} \label{ineq 1}
f\left(\frac{1}{2}\ub_0+\frac{1}{2}\ub_0'\right) < \frac{1}{2}f(\ub_0)+\frac{1}{2}f(\ub_0').
\end{equation}
If $\ub_0$ is a vertex (resp. a barycenter) of a cuboid $F + \lambda'$, then $\ub_0'$ is also of the same kind, therefore $\frac{1}{2}\ub_0+\frac{1}{2}\ub_0'$ is a vertex of $\Ccal$ and \eqref{ineq 1} follows from strict convexity of $q$. Otherwise, we have
\begin{equation} \label{typi 2}
\ub_0=\frac{1}{2}\left(b_1'+\dots+b_i'\right), \quad \ub_0'=\frac{1}{2}\left(b_1'+ \dots b_{i-1}'+b_{i+1}'\right) \quad (i \in \{2, \dots,n\})
\end{equation}
up to the operations (I)--(III). Then $\frac{1}{2}(\ub_0+\ub_0')=\frac{1}{2}(\ub_1+\ub_1')$ for the vertices $\ub_1=\frac{1}{2}(b_1'+\dots +b_{i-1}')$ and $\ub_1'=\frac{1}{2}(b_1'+\dots +b_{i+1}')$ of the same maximal simplex of $\Ccal$, hence 
\begin{equation} \label{eq 3}
f\left(\frac{1}{2}\ub_0+\frac{1}{2}\ub_0'\right)=\frac{1}{2}f(\ub_1)+\frac{1}{2}f(\ub_1')= \frac{1}{2}f(\ub_0) + \frac{1}{2} f(\ub_0'),
\end{equation}
where we have used in the last step that $f$ agrees with $q+\ell$ on the vertices of $\Ccal$ and that the vectors $b_1',\dots,b_n'$ are orthogonal  with respect to $b$. We conclude that \eqref{ineq 1} is {\it not} satisfied for $f$. 

We choose a new $f$ which is affine on the maximal simplices $\Delta$ of $\Ccal$ and which is a small perturbation of the old $f$. Let $\ve$ be a small positive rational number. We have to fix the values of $f$ on the vertices of the triangulation. For the simplex $\Delta$ in \eqref{barycentric}, let $f({\mathbf 0})=(q+\ell)({\mathbf 0})$ and let
$$f\left(\frac{1}{2}(b_1'+\dots+b_i')\right)=(q+\ell)\left(\frac{1}{2}(b_1'+\dots+b_i')\right)+\ve(2^{-1}+ \dots + 2^{-i})$$
for $i \in \{1,\dots,n\}$. Similarly, we define $f$ on the other maximal simplices with perturbation according to the distance to $\Lambda'$. Then $f$ still satisfies \eqref{periodicity}. We have to prove \eqref{ineq 1}. If $\ve$ is small enough, then \eqref{ineq 1} still holds for $\ub_0$ a vertex (resp. the barycenter) of a maximal cuboid $F+\lambda'$, $\lambda' \in \Lambda'$. In the situation \eqref{typi 2}, we deduce from \eqref{eq 3}
$$2f\left(\frac{1}{2}\ub_0+\frac{1}{2}\ub_0'\right)=
(q+\ell)(\ub_0)+(q+\ell)(\ub_0')
+
\ve\left(2^{-0}+\dots+2^{-(i-2)}+2^{-i}+2^{-(i+1)}\right)$$
and
$$f(\ub_0) + f(\ub_0')=
(q+\ell)(\ub_0)+(q+\ell)(\ub_0')
+
\ve\left(2^{-0}+\dots+2^{-(i-1)}\right)$$
which proves \eqref{ineq 1}. We conclude that $f$ is a strongly polyhedral convex function with respect to $\Ccal$. \qed

\begin{art} \rm \label{setting tropical equidistribution}
Let $X$ be an irreducible closed subvariety of $A$ of dimension $d$. The main theorem of tropical analytic geometry says that $\valbar(X_v^{\rm an})$ is a rational polytopal set in $\rtor$ of pure dimension $d$ (\cite{Gu3}, Theorem 6.9). 

For $P \in X(\overline{K})$, the ${\rm Gal}(\overline{K}/K)$-orbit is denoted by $O(P)$. 
Let $(P_m)_{m \in I}$ be a net in $X(\overline{K})$, i.e. $I$ is a directed set  and $P_m \in X(\overline{K})$. It is called a {\it generic net} in $X$ if for every proper closed subset $Y$ of $X$, there is $m_0 \in I$ such that $P_m \not \in Y$ for all $m \geq m_0$. The Dirac measure in $\ub \in \rdop^n$ is denoted by $\delta_\ub$.

We choose an embedding $\overline{K} \hookrightarrow \kdop_v$ to identify $A(\overline{K})$ with a subset of $A_v^{\rm an}= ({\Tor})_v^{\rm an}/M$. 
Now we are ready to state the {\it tropical equidistribution theorem}:
\end{art} 

\begin{thm} \label{tropical equidistribution theorem}
Let $L$ be an ample symmetric line bundle on $A$ endowed with a canonical $M_B$-metric $\metr_{\rm can}$. Let $(P_m)_{m \in {I}}$ be a generic net in $X(\overline{K})$ with $\lim_{m } \hat{h}_L(P_m)=0$. Then we have
\begin{equation} \label{tropical equidistribution}
\frac{1}{|O(P_m)|} \cdot \sum_{P_m^\sigma \in O(P_m)} \delta_{\valbar(P_m^\sigma)} \stackrel{w}{\rightarrow} 
\mu:= \frac{1}{\deg_L(X)}\cdot\valbar\left(c_1(L|_X,\metr_{{\rm can},v})^{\wedge d}\right) 
\end{equation}
as a weak limit of regular probability measures on the tropical variety $\valbar(X_v^{\rm an})$. Moreover, the equidistribution measure is a strictly positive piecewise Haar measure on the polytopal set $\valbar(X_v^{\rm an})$.
\end{thm}

The last claim means that $\valbar(X_v^{\rm an})$ is a finite union of $d$-dimensional polytopes $\Deltabar$ such that $\mu|_\Deltabar$ is a positive multiple of the Lebesgue measure. This follows from the first claim and \cite{Gu3}, Theorem 1.3. 

The condition $\lim_{m} \hat{h}_L(P_m)=0$ is independent of the choice of the ample symmetric line bundle $L$ and hence the theorem shows that $\mu$ is completely independent of $L$.

\proof We have to test weak convergence for a continuous function $f$ on $\valbar(X_v^{\rm an})$. By weak compactness of the set of regular probability measures on $\valbar(\Xan)$ (see \cite{BG}, 4.3), it is enough to show that every convergent subnet of the left hand side of \eqref{tropical equidistribution} converges to the right hand side. Passing to a subnet, we may assume that we have weak convergence
\begin{equation} \label{weak subconvergence}
\frac{1}{|O(P_m)|} \cdot \sum_{P_m^\sigma \in O(P_m)} \delta_{\valbar(P_m^\sigma)} \stackrel{w}{\rightarrow}  \nu 
\end{equation}
for a regular probability measure $\nu$ on $\valbar(X_v^{\rm an})$ and we have to show $\mu=\nu$. Let $\Ccalbar$ be the triangulation of $\rtor$ from Lemma \ref{barycentric model} and let $\Ccal_j:= 2^{-j}\Ccal$ for $j \in \ndop$. We assume that our test function $f$ is a piecewise affine $\Lambda$-periodic function with
\begin{equation} \label{affine equation 1}
f(\ub)=\mb_\Delta \cdot \ub + c_\Delta
\end{equation}
on every $n$-dimensional $\Delta \in \Ccal_j$ for suitable $\mb_\Delta \in \zdop^n$, $c_\Delta \in \qdop$. Note that such functions, with varying $j \in \ndop$, generate a dense $\qdop$-subspace of $C(\valbar(\Xan))$. By \eqref{weak subconvergence}, it is enough to prove
\begin{equation} \label{measure equation}
\int_{\valbar(X_v^{\rm an})} f\,d\mu = \int_{\valbar(X_v^{\rm an})} f\,d\nu.
\end{equation}
Now $j$ is fixed for our given $f$. By Mumford's construction, there is a formal $\kdop_v^\circ$-model $\Acal_v^{(j)}$ of $A$ associated to the triangulation $\Ccalbar_j$ of $\rtor$. For $g:= f \circ \valbar$, the metric $\metr_{g,v}$ from \ref{Chern forms} is a formal metric with respect to a formal $\kdop_v^\circ$-model $\Gcal_v$ of $O_X$ on $\Acal_v^{(j)}$ (\cite{Gu3}, Proposition 6.6). 
Let $\Lcal_v$ be the $\kdop_v^\circ$-model of $L^{\otimes m_j}$ on $\Acal_v^{(j)}$ considered in Lemma \ref{barycentric model}. Replacing $L$ by $L^{\otimes m_j}$, we may assume $m_j=1$. The associated formal metric on $L$ is denoted by $\metr_0$. As in \ref{toric line bundles}, we identify $p^*(L_v^{\rm an})$ with the trivial line bundle on $T$. There is a unique rational piecewise affine continuous function $f_0:\rdop^n \rightarrow \rdop$ with
\begin{equation} \label{f characterization}
f_0 \circ \val=-\log(p^*||1||_0)
\end{equation} 
on $T$ (\cite{Gu3}, Proposition 6.6). Ampleness of $\tilde{\Lcal}_v$ corresponds to the fact that $f_0$ is a strongly polyhedral convex function with respect to $\Ccal_j$ (see \cite{Gu3}, Corollary 6.7). We have
\begin{equation} \label{affine equation for f_0}
f_0(\ub) = \mb_\Delta^{(0)} \cdot \ub + c_\Delta^{(0)}
\end{equation}
on the $n$-dimensional elements $\Delta$ of $\Ccal_j$ for suitable $\mb_\Delta^{(0)} \in \zdop^n$, $c_\Delta^{(0)} \in \qdop$. Note that a piecewise affine function $f_0$ with \eqref{affine equation for f_0} is a strongly polyhedral convex function if and only if, for $n$-dimensional simplices $\Delta, \sigma \in \Ccal_j$ with $\codim(\Delta \cap \sigma, \rdop^n)=1$, we have
\begin{equation} \label{strong inequality}
\nb_{\Delta,\sigma} \cdot \mb_\Delta^{(0)} > \nb_{\Delta,\sigma} \cdot \mb_\sigma^{(0)}
\end{equation}
where $\nb_{\Delta, \sigma}$ is the inner normal vector of $\Delta$ at $\Delta \cap \sigma$ of length $1$. By \eqref{periodicity}, \eqref{cycle and b} and \eqref{affine equation for f_0}, the quantity $\nb_{\Delta,\sigma} \cdot (\mb_\Delta^{(0)}-\mb_\sigma^{(0)})$ remains invariant if we replace $(\Delta,\sigma)$ by $(\Delta + \lambda, \sigma + \lambda)$ for $\lambda \in \Lambda$. Hence we may choose a sufficiently small positive rational number $q$ with
\begin{equation} \label{q inequality}
q \max_{(\Delta,\sigma)} |\nb_{\Delta,\sigma}\cdot(\mb_\Delta-\mb_\sigma)| <
\min_{(\Delta,\sigma)} \nb_{\Delta,\sigma} \cdot (\mb_\Delta^{(0)}-\mb_\sigma^{(0)}),
\end{equation}
where $(\Delta,\sigma)$ ranges over all pairs as above. Therefore \eqref{strong inequality} holds for $f_0+qf$ and we conclude that $f_0+qf$ is a rational strongly polyhedral convex function with respect to $\Ccal_j$. This means that some positive tensor power of the metric
$$\metr_{0,v}':= \metr_{0,v} \otimes \metr_{qg,v}$$
is a formal metric associated to a formal $\kdop_v^\circ$-model on $\Acal_v^{(j)}$ with ample reduction. We extend the metrics $\metr_{0,v}$ and $\metr_{0,v}'$ to admissible semipositive $M_B$-metrics on $L$ by using the given canonical metrics $\metr_{{\rm can},w}$ at the other places $w$ of $M_B$. 

We choose a rigidification of $L$ associated to the canonical metric $\metr_{\rm can}$. By the theorem of the cube, we may identify $[m]^*(L)$ with $L^{\otimes m^2}$. By the  construction of canonical metrics, we have 
$$\metr_{\rm can} = \lim_{i \to \infty} \metr_i, \quad \metr_i:=\left([2^i]^* \metr_0 \right)^{4^{-i}}.$$
Note that the limit occurs only at the fixed place $v$. More precisely, we have
\begin{equation} \label{metric distance}
d(\metr_{i,v}, \metr_{{\rm can},v}) = O(4^{-i})
\end{equation} 
with respect to the distance from \ref{isometry classes}. This follows from the fact that $\metr_v \mapsto ([2]^* \metr_v)^{1/4}$ is a contraction with factor $\frac{1}{4}$ on the space of bounded metrics on $L_v^{\rm an}$ (see \cite{BG}, proof of Theorem 9.5.4). In \eqref{metric distance} and in the following estimates, the implicit constant is always independent of $i$ and of the net $(P_m)_{m \in {I}}$,  but may depend on the geometric data and $f$. 

We repeat now the above considerations with the metric $\metr_{i,v}$ replacing $\metr_{0,v}$ (i.e. $j$ is replaced by $i+j$). The morphism $[2^i]$ extends to a morphism $\varphi_i:\Acal_v^{(i+j)} \rightarrow \Acal_v^{(j)}$ over $\kdop_v^\circ$ with finite reduction (see \cite{Gu3}, Proposition 6.4). Then $\varphi_i^*(\Lcal_v)$ is a $\kdop_v^\circ$-model of $[2^i]^*(L)=L^{\otimes 4^i}$ on $\Acal_v^{i+j}$ with formal metric equal to $\metr_{i,v}^{\otimes 4^i}$. We conclude that the function $f_i$ associated to $\metr_i$ by \eqref{f characterization} (with $i$ replacing $0$) is a strongly polyhedral convex function with respect to $\Ccal_{i+j}$ satisfying 
\begin{equation} \label{f_i properties}
f_i(\ub)=4^{-i}f_0(2^i\ub), \quad \mb_{\Delta'}^{(i)}=2^{-i}\mb_{2^i\Delta'}^{(0)}
\end{equation}
for every $n$-dimensional simplex $\Delta'$ of $\Ccal_{i+j}$. We claim that $f_i+2^{-i}qf$ is a rational strongly polyhedral convex function with respect to $\Ccal_{i+j}$. To show this, let $\Delta',\sigma'$ be $n$-dimensional simplices of $\Ccal_{i+j}$ such that $\Delta' \cap \sigma'$ is of codimension $1$ in both. Since $\Ccal_{i+j}$ is a refinement of $\Ccal_j$, there are $\Delta, \sigma \in \Ccal_j$ with $\Delta' \subset \Delta$, $\sigma' \subset \sigma$. We have to show that
\begin{equation} \label{strong inequality 2}
\nb_{\Delta',\sigma'} \cdot \left(\mb_{\Delta'}^{(i)}+2^{-i}q\mb_\Delta-\mb_{\sigma'}^{(i)}-2^{-i}q\mb_\sigma \right) > 0.
\end{equation}
If $\Delta=\sigma$, then this follows from \eqref{f_i properties} and \eqref{strong inequality}. If $\Delta \neq \sigma$, then $\Delta \cap \sigma$ is of codimension $1$ in $\rdop^n$ and $\nb_{\Delta',\sigma'}=\nb_{\Delta,\sigma}$. We conclude that \eqref{strong inequality 2} follows from \eqref{f_i properties} and \eqref{q inequality}. Hence $f_i$ is a strongly polyhedral convex function with respect to $\Ccal_{i+j}$. 

We use this to define an admissible semipositive $M_B$-metric $\metr_i'$ on $L$ by 
$$\metr_{i,v}':= \metr_{i,v} \otimes \metr_{2^{-i}qg,v}$$
and using $\metr_{{\rm can},w}$ at the other places. By Theorem \ref{properties of local height}, we have
\begin{equation} \label{equi 1}
h_{(L, \metr_i')}(X)= h_{(L, \metr_i)}(X) + \sum_{r=1}^{d+1} \binom{d+1}{r}(2^{-i}q)^r h_r^{(i)},
\end{equation}
where
$$h_r^{(i)}:= \lambda_{\underbrace{(O_X,\metr_{g,v}). \dots, (O_X, \metr_{g,v})}_r, \underbrace{(L,\metr_{i,v}),\dots, (L, \metr_{i,v})}_{d+1-r}}(X,v)$$
for $r \in \ndop$. If we replace the metric $\metr_{i,v}$ by $\metr_{{\rm can},v}$ in the definition of $h_r^{(i)}$, then we get a quantity $\hat{h}_r^{(i)}$. Note that the formal metric $\metr_{g,v}$ is the quotient of two semipositive formal metrics on an ample line bundle. By Theorem \ref{properties of local height} and \eqref{metric distance}, we get
$$h_r^{(i)}= \hat{h}_r^{(i)} + O(4^{-i}).$$
Using this in \eqref{equi 1} for $r=0, \dots, d+1$ and $\hat{h}_L(X)=0$ (Lemma \ref{fundamental inequality}), we get
\begin{equation} \label{equi 2}
h_{(L, \metr_i')}(X)= (d+1)2^{-i}q \int_{\Xan} g \cdot c_1(L|_X,\metr_{{\rm can},v})^{\wedge d} + O(4^{-i}).
\end{equation}
The transformation formula and the definitions of $\mu,g$ imply
\begin{equation} \label{equi 3}
\int_{\Xan} g \cdot c_1(L|_X,\metr_{{\rm can},v})^{\wedge d}= \deg_L(X) \int_{\valbar(X_v^{\rm an})} f \, d\mu.
\end{equation}
Since $\metr_i'$ is semipositive, the fundamental inequality (Lemma \ref{fundamental inequality}) yields 
\begin{equation} \label{equi 4}
\frac{1}{(d+1)\deg_L(X)} h_{(L,\metr_i')}(X) \leq \liminf_{m } h_{(L,\metr_i')}(P_m).
\end{equation}
By linearity and \eqref{metric distance}, we have
\begin{equation} \label{equi 5}
h_{(L,\metr_i')}(P_m)= \hat{h}_L(P_m)+ \frac{q}{2^i|O(P_m)|} \sum_{P_m^\sigma \in O(P_m)} g(P_m^\sigma) + O(4^{-i}).
\end{equation}
Using $\hat{h}_L(P_m) \to 0$, \eqref{weak subconvergence} and \eqref{equi 4} in \eqref{equi 5}, we deduce 
\begin{equation} \label{equi 6}
\frac{1}{(d+1)\deg_L(X)} h_{(L,\metr_i')}(X) \leq 2^{-i}q \int_{\valbar(X_v^{\rm an})} f \, d\nu + O(4^{-i}).
\end{equation}
Finally, we put \eqref{equi 3} and \eqref{equi 6} in \eqref{equi 2} to get
$$2^{-i}q(d+1) \deg_L(X) \int_{\valbar(X_v^{\rm an})} f \, d\nu \geq 2^{-i} q (d+1) \deg_L(X) \int_{\valbar(X_v^{\rm an})} f \, d\mu + O(4^{-i}).$$
For $i \to \infty$, this is only possible if ``$\leq$'' holds in \eqref{measure equation}. Replacing $f$ by $-f$, we get ``$=$'' and hence $\mu=\nu$. \qed

\begin{rem} \label{Moriwaki} \rm
The results of this section hold more generally for a fixed discrete valuation $v$ of a field with product formula for which the fundamental inequality holds. In particular, this is true for number fields (\cite{Zh1}, Theorem 5.2) and more generally for the finitely generated fields over $\qdop$ with the product formula considered by Moriwaki (see \cite{Mo4}), Corollary 5.2). 
\end{rem}

\begin{rem} \label{other proofs} \rm 
There are similar equidistribution theorems by Szpiro, Ullmo and Zhang (\cite{SUZ}, Th\'eor\`eme 3.1) in the number field case at an archimedean place $v$, by Moriwaki (\cite{Mo4}, Theorem 6.1) in the generalization mentioned in Remark \ref{Moriwaki} ($v$ again archimedean) and by Chambert-Loir (\cite{Ch}, Th\'eor\`eme 3.1) in the number field case at a finite $v$ with respect to an ample formal metric on $L$ at $v$. Their proofs were less subtle than the argument for Theorem \ref{tropical equidistribution theorem} as  the twist $\metr_{i,v}'=\metr_{i,v} \otimes \metr_{qg,v}$ was already semipositive for $q$ sufficiently small independent of $i$. Therefore the uniformity \eqref{metric distance} and the careful choice of the models in Lemma \ref{barycentric model} played no role. 
\end{rem}

\section{The Bogomolov conjecture}

Let $K=k(B)$ be a function field. We prove the Bogomolov conjecture for an abelian variety $A$ over $K$ which is totally degenerate at some place $v \in M_B$. The proof follows closely Zhang's proof in the number field case but the dimensionality argument is now using the associated tropical varieties. At the end, we give some applications  analogous to the corollaries in Zhang's paper \cite{Zh2}. We start with a lemma which holds in more generality:


\begin{lem} \rm \label{totally degenerate properties}
Let $F$ be a field with a non-trivial non-archimedean complete absolute value $v$ and let $A$ be a totally degenerate abelian variety over $F$.
\begin{itemize}
\item[(a)] Every abelian subvariety of $A$ is totally degenerate.
\item[(b)] If $\varphi:A \rightarrow B$ is a surjective homomorphism of abelian varieties over $F$, then $B$ is also totally degenerate.
\item[(c)] If $\varphi:C\rightarrow A$ is a homomorphism of abelian varieties over $F$ and if $C$ has good reduction, then $\varphi \equiv 0$.
\end{itemize}
\end{lem}

\proof If $\varphi: A_1 \rightarrow A_2$ is a homomorphism of abelian varieties over $F$, then the associated Raynaud extensions are homomorphic (see \cite{BL1}, \S 1 and use the argument before Proposition 3.5). If $\varphi$ is an isogeny, then it is clear that the corresponding tori (resp. abelian varieties) of the Raynaud extensions are also isogeneous. By Poincar\'e's complete reducibility theorem (\cite{Mu2}, IV.19, Theorem 1), ${\rm ker}(\varphi) \times \varphi(A_1)$ is isogeneous to $A_1$. This proves easily (a)--(c). \qed

\vspace{3mm}

\label{proof of Bogo} \noindent {\bf Proof of Theorem \ref{Bogomolov conjecture}: \/} 
By Lemma \ref{totally degenerate properties}(c),  the Chow trace ${\rm Tr}_{K'/k}(A)$ of $A$ is $0$ for every finite extension $K'/K$ and hence the points of N\'eron--Tate height $0$ are the torsion points in $A(\overline{K})$ (see \cite{La}, Theorem 6.5.4). This proves the claim in the $0$-dimensional case.

For $d:=\dim(X)>0$, we argue by contradiction. We may assume that the counterexample $X$ is irreducible. Let $G$ be the stabilizer of $X$ in $A$. Then $B:=A/G$ is an abelian variety. Let $\pi:A \rightarrow B$ be the quotient homomorphism. The closed subvariety $X':=\pi(X)$ of $B$ has trivial stabilizer. Since $B$ is totally degenerate at $v$ (Lemma \ref{totally degenerate properties}(b)), $X'$ is also a counterexample to the Bogomolov conjecture. Indeed, this follows from $\hat{h}_{\pi^*(L')} \ll \hat{h}_L$ for every ample symmetric line bundle $L'$ on $B$. Hence we may assume that the stabilizer of $X$ is trivial. By finite base extension, we may assume that $X$ is defined over $K$. 
For $N \in \ndop$ sufficiently large, the argument in \cite{Ab}, Lemma 4.1, shows that the morphism
$$\alpha:X^N \rightarrow A^{N-1}, \quad (x_1, \dots , x_N) \mapsto (x_2-x_1, \dots, x_N- x_{N-1})$$
is generically finite. 
On $A^N$, we use the N\'eron--Tate height $\hat{h}$ with respect to the ample symmetric line bundle $L_N:=p_1^*(L) \otimes \cdots \otimes p_N^*(L)$. 
We have 
$$X(\ve)^N \subset X^N(N\ve)$$
for all $\ve >0$ and therefore $X^N$ is also a counterexample for the Bogomolov conjecture. Hence there is a generic net $(\xb_m)_{m \in I}$ in $\left(X^N\right)(\overline{K})$ with $\lim_m \hat{h}(\xb_m)=0$. Obviously, $(\alpha(\xb_m))_{m \in I}$ is also a generic net in $Y:=\alpha(X^N)$ with N\'eron--Tate height converging to $0$. 
For both nets, we may apply the tropical equidistribution theorem to get strictly positive Haar measures $\mu$, $\nu$ on $\valbar(\Xan)^N$ and $\valbar(Y_v^{\rm an})$ with
\begin{equation} \label{Bogo1}
\lim_{m } \frac{1}{|O(\xb_m)|} \sum_{\xb_m^\sigma \in O(\xb_m)} \delta_{\valbar(\xb_m^\sigma)} \stackrel{w}{\rightarrow} \mu, \quad
\lim_{m } \frac{1}{|O(\alpha\xb_m)|} \sum_{\yb_m^\sigma \in O(\alpha\xb_m)} \delta_{\valbar(\yb_m^\sigma)} \stackrel{w}{\rightarrow} \nu.
\end{equation}

We have the linear map
$$\alpha_\val: (\rdop^n)^N \longrightarrow (\rdop^n)^{N-1}, \quad (\ub_1, \dots, \ub_N) \mapsto (\ub_2-\ub_1, \dots, \ub_N-\ub_{N-1})$$
which satisfies
\begin{equation} \label{Bogo2}
\overline{\alpha}_\val \circ \valbar = \valbar \circ \alpha
\end{equation}
on $(A_v^{\rm an})^N$. Clearly, we have
$$(\overline{\alpha}_\val)_* \left(\frac{1}{|O(\xb_m)|} \sum_{\xb_m^\sigma \in O(\xb_m)} \delta_{\valbar(\xb_m^\sigma)} \right)
= \frac{1}{|O(\alpha\xb_m)|} \sum_{\yb_m^\sigma \in O(\alpha\xb_m)} \delta_{\valbar(\yb_m^\sigma)}$$
and hence $(\overline{\alpha}_\val)_*(\mu)=\nu$ by \eqref{Bogo1}. We consider $X$ as a closed subvariety of $X^N$ by using the diagonal map. Note that $\alpha(X)=\{0\}$ in $A^{N-1}$. By \cite{Gu3}, Theorem 6.9, $\valbar(\Xan)$ is a polytopal set in $\rdop^{Nn}$ of pure dimension $d$. Moreover, there is an $Nd$-dimensional simplex $\Deltabar$ in $\valbar(\Xan)^N$ with a $d$-dimensional face contained in the diagonal  $\valbar(\Xan)$. By \eqref{Bogo2}, $\overline{\alpha}_\val$ maps this face to $\mathbf 0$ and hence $\overline{\tau} = \overline{\alpha}_\val(\Deltabar)$ is a simplex of dimension $\leq(N-1)d$. Since $\nu$ is a Haar measure on $\valbar(Y_v^{\rm an})$, we have $\nu(\overline{\tau})=0$. For every $\ve >0$, there is a continuous function $f$ on $\valbar(Y_v^{\rm an})$ with $0 \leq f \leq 1$, $f({\mathbf 0})=1$ and $\int f \, d\nu < \ve$. We conclude that 
$$\mu(\Deltabar) \leq \int_{\valbar(\Xan)^N} f \circ \overline{\alpha}_\val \, d\mu = \int_{\valbar(Y_v^{\rm an})} f \, d\nu < \ve$$
and hence $\mu(\Deltabar)=0$. This contradicts strict positivity of the Haar measure $\mu$. \qed 

\begin{art} \label{torsion subvarieties} \rm 
A {\it torsion subvariety} of $A$ is a translate of an abelian subvariety of $A$ by a torsion point. Again, we consider a closed subvariety $X$ of $A$ defined over $\overline{K}$. Let $X^*$ be the complement in $X$ of the union of all torsion subvarieties contained in $X$.
\end{art}

\begin{thm} \label{Bogomolov conjecture 2} 
The closed subset $X \setminus X^*$ of $X$ is a finite union of torsion subvarieties which are maximal in $X$ and we have
$$\inf\{\hat{h}_L(P)\mid P \in X^*(\overline{K})\}>0$$
for any ample symmetric line bundle $L$ on $A$.
\end{thm}

\proof This follows from Theorem \ref{Bogomolov conjecture} as in the number field case (see \cite{Ab}, Th\'eor\`eme 4.3). \qed

\begin{art} \rm \label{rationality}
If $X$ is defined over a finite extension $F/K$, then the conjugates over $F$ of a maximal torsion subvariety contained in $X$ are of the same kind and hence $X^*$ is defined over $F$ (see \cite{BG}, A.4.13). 
\end{art}

\begin{rem} \label{applications} \rm 
For the following applications, which were suggested by the referee, we denote by $K$ a function field as before or a number field (resp. more generally a finitely generated field over $\qdop$ with the product formula relative to a big polarization introduced by Moriwaki \cite{Mo4}). By Remark \ref{Moriwaki}, the above proof of Bogomolov's conjecture  holds also for the latter case. This provides only a non-archimedean proof of a special case of known results (see Section 1). However, the non-archimedean method leads to the following applications:

We still assume that the abelian variety $A$ over  $K$ is to\-tal\-ly degenerate at the non-archimedean place $v$. We choose an embedding $\overline{K} \hookrightarrow \kdop_v$ to identify $A(\overline{K})$ with a subset of $A_v^{\rm an}= ({\Tor})_v^{\rm an}/M$ and we will use the map $\valbar:A_v^{\rm an} \rightarrow \rtor$ from Section 5. We choose also an ample symmetric line bundle $L$ on $A$.
\end{rem}


\begin{cor} \label{strict small sequences}
Let $(P_m)_{m \in \ndop}$ be a sequence in $A(\overline{K})$ such that no infinite subsequence is contained in a proper torsion subvariety of $A$. If $\lim_{m \to \infty} \hat{h}_L(P_m)=0$, then the associated ${\rm Gal}(\overline{K}/K)$-orbits $O(P_m)$ are tropically equidistributed on $\valbar(A_v^{\rm an})=\rtor$, i.e.
$$\frac{1}{|O(P_m)|} \cdot \sum_{P_m^\sigma \in O(P_m)} \delta_{\valbar(P_m^\sigma)} \stackrel{w}{\rightarrow} 
\mu$$
for the Haar probability measure $\mu$ on $\rtor$.
\end{cor}

\proof It follows from \cite{Gu3}, Corollary 9.9, that $\deg_L(A)\mu=\valbar(c_1(L,\metr_{\rm can})^{\wedge n})$. By Theorem \ref{tropical equidistribution theorem}, it is enough to prove that $(P_m)_{m \in \ndop}$ is a generic sequence in $A$. We argue by contradiction and assume that there is an infinite subsequence $(P_m)_{m \in I}$ dense in a proper subvariety $X$ of $A$. By assumption, $X$ is no torsion subvariety of $A$. By passing to a subsequence, we may assume that  $P_m \in X^*$ for all $m \in I$. We have still $\lim_{m \to \infty} \hat{h}_L(P_m)=0$ contradicting Theorem \ref{Bogomolov conjecture 2}. \qed



\begin{cor} \label{totally split points}
Let $K^{nr}$ be the maximal subextension of $\overline{K}$ which is unramified over $v$ and let ${K'}$ be a finite extension of $K^{nr}$.  Then the number of torsion points in $A({K'})$ is finite and there is $\ve >0$ such that $\hat{h}_L(P) \geq \ve$ for all non-torsion points $P \in A({K'})$. 
\end{cor}

\proof We prove the claim by contradiction. We may assume that $A$ has minimal dimension such that the result is wrong. 
We conclude that $A({K'})$ contains infinitely many different points $(P_m)_{m \in \ndop}$ with $\lim_{m \to \infty} \hat{h}_L(P_m) =0$. We may assume that ${K'}$ is a normal subfield of $\kdop_v$.

We may suppose that the uniformization $A_v^{\rm an} = ({\Tor})_v^{\rm an}/M$ is analytically defined over a subfield $F$ of $\kdop_v$ such that $F$ is finite over the completion $K_v$ of $K$ (see \cite{BL1}, Section 1). 
Then we have $A(E)=(E^\times)^n/M$ for every algebraic subextension $E/F$ of $\kdop_v$. 
Since  $K^{nr}F$ is unramified over $F$, we conclude that the ramification index $e$ of ${K'}F$ over $K_v$ is finite. We note that $\val(A({K'}))$ is contained in $\frac{1}{e}\zdop^n$. 

It follows that the sequence $\valbar(O(P_m))$ can't be equidistributed in $\rtor$. By Corollary \ref{strict small sequences}, there is an infinite subsequence $(P_m)_{m \in I}$ contained in $\zeta + B$ for a torsion point $\zeta$ and a proper abelian subvariety $B$.  For the finite extension $K'':=K'(\zeta)$ of $K^{nr}$, the infinite sequence $Q_m:= P_m - \zeta$, $m \in I$, is contained in $B(K'')$. Since $\hat{h}_L(Q_m)=\hat{h}_L(P_m)$ tends to $0$, we conclude that the corollary can't be true for $B$. This contradicts the minimality of $A$. \qed


{\small Walter Gubler, Fachbereich Mathematik, Universit\"at Dortmund,
 D-44221 Dortmund, walter.gubler@mathematik.uni-dortmund.de}
\end{document}